
%
\magnification=\magstephalf
\documentstyle{amsppt}
 \addto\tenpoint{\baselineskip 15pt
  \abovedisplayskip18pt plus4.5pt minus9pt
  \belowdisplayskip\abovedisplayskip
  \abovedisplayshortskip0pt plus4.5pt
  \belowdisplayshortskip10.5pt plus4.5pt minus6pt}\tenpoint
\pagewidth{6.5truein} \pageheight{8.9truein}
\subheadskip\bigskipamount
\belowheadskip\bigskipamount
\aboveheadskip=3\bigskipamount
\catcode`\@=11
\def\output@{\shipout\vbox{%
 \ifrunheads@ \makeheadline \pagebody
       \else \pagebody \fi \makefootline
 }%
 \advancepageno \ifnum\outputpenalty>-\@MM\else\dosupereject\fi}
\outer\def\subhead#1\endsubhead{\par\penaltyandskip@{-100}\subheadskip
  \noindent{\subheadfont@\ignorespaces#1\unskip\endgraf}\removelastskip
  \nobreak\medskip\noindent}
\outer\def\enddocument{\par
  \add@missing\endRefs
  \add@missing\endroster \add@missing\endproclaim
  \add@missing\enddefinition
  \add@missing\enddemo \add@missing\endremark \add@missing\endexample
 \ifmonograph@ 
 \else
 \vfill
 \nobreak
 \thetranslator@
 \count@\z@ \loop\ifnum\count@<\addresscount@\advance\count@\@ne
 \csname address\number\count@\endcsname
 \csname email\number\count@\endcsname
 \repeat
\fi
 \supereject\end}
\catcode`\@=\active
\CenteredTagsOnSplits
\NoBlackBoxes
\nologo
\def\today{\ifcase\month\or
 January\or February\or March\or April\or May\or June\or
 July\or August\or September\or October\or November\or December\fi
 \space\number\day, \number\year}
\define\({\left(}
\define\){\right)}
\define\Ahat{{\hat A}}

\define\CC{{\Bbb C}}

\define\Hom{\operatorname{Hom}}

\define\RR{{\Bbb R}}

\define\Tr{\operatorname{Tr}}
\define\ZZ{{\Bbb Z}}
\define\[{\left[}
\define\]{\right]}

\define\chiup{\raise.5ex\hbox{$\chi$}}
\define\cir{S^1}

\define\exertag #1#2{#2\ #1}

\define\index{\operatorname{index}}

\define\inv{^{-1}}
\define\mstrut{^{\vphantom{1*\prime y}}}
\define\protag#1 #2{#2\ #1}

\define\res#1{\negmedspace\bigm|_{#1}}
\define\temsquare{\raise3.5pt\hbox{\boxed{ }}}

\define\theprotag#1 #2{#2~#1}

\define\xca#1{\removelastskip\medskip\noindent{\smc%
#1\unskip.}\enspace\ignorespaces }

\define\zmod#1{\ZZ/#1\ZZ}

\redefine\Re{\operatorname{Re}}

\NoRunningHeads
\define\Det{\operatorname{Det}}
\define\Ker{\operatorname{Ker}}
\define\Rnz{\RR^{\not= 0}}
\define\YG{Y_\Gamma }
\define\Ygh{\hat{Y}_\gamma }
\define\Yg{Y_\gamma }
\define\alim{\operatorname{a-lim}}
\define\bX{\partial X}

\define\comp#1#2{\left[ #1 \right]_{(#2)}}
\define\curv#1{\Omega^{#1}}
\define\cut{^{\text{cut}}}

\define\dibf{\Det\inv _{\bX/Z}}
\define\dib{\Det\inv _{\bX}}
\define\hol{\operatorname{hol}}

\define\sign{\operatorname{sign}}
\define\sind{(-1)^{\index D_Y}}
\define\spec{\operatorname{spec}}
\define\spin{\operatorname{spin}}
\define\tpi{2\pi i}
\define\zo{[0,1]}

\input epsf

\refstyle{A}
\widestnumber\key{SSSSSSSS}   
\document

	\topmatter
 \title\nofrills Determinant Line Bundles Revisited \endtitle
 \author Daniel S. Freed  \endauthor
 \thanks To appear in the proceedings of the conference {\it Topological and
Geometrical Problems related to Quantum Field Theory,\/} Trieste, Italy,
March 13--24, 1995.\endgraf
 The author is supported by NSF grant DMS-9307446, a Presidential
Young Investigators award DMS-9057144, and by the O'Donnell
Foundation.\endthanks
 \affil Department of Mathematics \\ University of Texas at Austin\endaffil
 \address Department of Mathematics, University of Texas, Austin, TX
78712\endaddress
 \email dafr\@math.utexas.edu \endemail
 \date May 11, 1995\enddate
	\endtopmatter

\document

This note is an addendum to joint work with Xianzhe
Dai~\cite{DF1},~\cite{DF2}.\footnote{In~\cite{DF1} the reader will find an
extensive discussion of related work and a bibliography.} In that paper we
investigate the geometric theory of {\it $\eta $-invariants\/} of Dirac
operators on manifolds with boundary.  We summarize the main results below.
One key geometric observation is that the exponentiated $\eta $-invariant
naturally takes values in the {\it determinant line\/} of the boundary.  As
such it is intimately related to the geometry of determinant line bundles for
{\it families\/} of Dirac operators.  The differential geometry of
determinant line bundles was developed first by Quillen~\cite{Q} in a special
case, and then by Bismut and Freed~\cite{BF1}, ~\cite{BF2} in general.
(See~\cite{F1} for an exposition of these results.)  In~\S{5} of~\cite{DF1}
the results on $\eta $-invariants are used to reprove the holonomy formula
for determinant line bundles, also known as Witten's global anomaly
formula~\cite{W}.  However, the argument there is unnecessarily complicated.
The main purpose of this note, then, is to reprove {\it both\/} the curvature
and holonomy formulas for determinant line bundles using the results
of~\cite{DF1}.  (The argument was sketched in~\cite{DF2}.)

To avoid repetitious recitation of requirements, we set some conventions here
which apply throughout.  We work with {\it compact\/} Riemannian manifolds.
If the boundary is nonempty we assume that the metric is a product near the
boundary.  Our results hold for any Dirac operator on a $\spin^c$ manifold
coupled to a vector bundle with connection, but for simplicity we state the
formulas only for the basic Dirac operator on a spin manifold.  Thus all
manifolds are assumed spin.  We use the $L^2$~metric on the spinor
fields~$S$.  A {\it family of Riemannian manifolds\/} is a smooth fiber
bundle $\pi \:X\to Z$ together with a metric on the relative (vertical)
tangent bundle~$T(X/Z)$ and a distribution of ``horizontal'' complements
to~$T(X/Z)$ in~$TX$.  We assume that $T(X/Z)$~is endowed with a spin
structure.  Also, when working with families of manifolds with boundary, we
assume that the Riemannian metrics on the fibers are products near the
boundary.  There is an induced family $\partial \pi \:\partial X\to Z$ of
closed manifolds.  Finally, we will always use~`$X$' to denote an odd
dimensional manifold and `$Y$' to denote an even dimensional manifold.

\nobreak
As stated earlier this is a continuation of joint work with Xianzhe Dai.

\comment
lasteqno @ 26
\endcomment

 \subhead Eta Invariants on Manifolds with Boundary
 \endsubhead

First recall that on a {\it closed\/} odd dimensional manifold~$X$ the Dirac
operator~$D_X$ is self-adjoint and has a discrete spectrum~$\spec(D_X)$
extending to~$+\infty $ and~$-\infty$.  The $\eta $-invariant of
Atiyah-Patodi-Singer~\cite{APS} is defined by meromorphic continuation of
the function
  $$ \eta _X(s) = \sum\limits_{{\lambda \not= 0}\atop{\lambda \in
     \spec(D_X)}}\frac{\sign\lambda }{|\lambda |^s},  $$
which by general estimates converges for $\Re(s)$~sufficiently large.  In
fact, for Dirac operators the meromorphic continuation is analytic
for~$\Re(s)>-2$~\cite{BF2, Theorem~2.6}.  In any case $\eta _X$~is regular
at~$s=0$, and we set
  $$ \tau \mstrut _X = \exp \pi i(\eta _X(0) + \dim\Ker D_X) \in
     \CC.  \tag{1} $$
The general theory of $\eta $-invariants shows that $\tau \mstrut _X$~varies
smoothly in families, whereas the $\eta $-invariant~$\eta _X(0)$ is
discontinuous in general.  Note that $|\tau \mstrut _X|=1$.

On a manifold with boundary we need to specify elliptic boundary conditions
to obtain an operator with discrete spectrum.  We use the boundary conditions
introduced by Atiyah-Patodi-Singer, but adapted to odd dimensional
manifolds~$X$.  This involves an additional piece of information
concerning~$\Ker D_{\bX}$.  Recall that on an even dimensional manifold~$Y$
the spinor fields~$S\mstrut _Y$ split as $S\mstrut _Y=S_Y^+\oplus S_Y^-$, and
the Dirac operator $D\mstrut _Y\:S^\pm_Y\to S^\mp_Y$ interchanges the
positive and negative pieces.  (In the sequel we use `$D\mstrut _Y$'~to
denote the operator $D\mstrut _Y\:S_Y^+\to S_Y^-$.)  If $Y=\bX$ is the
boundary of an odd dimensional manifold~$X$, then $\dim\Ker^+D_{\bX} =
\dim\Ker^- D_{\bX}$.  The additional piece of information we must choose as
part of the boundary condition is an isometry
  $$ T\:\Ker^+D_{\bX}\longrightarrow \Ker^-D_{\bX}.  $$
Then the basic analytic properties of~$D_X$ with these boundary conditions
are the same as those of the Dirac operator on a closed manifold, and so the
invariant~\thetag{1} is defined. Its dependence on~$T$ is simple, and
factoring this out we observe that
  $$ \tau \mstrut _X\in \dib, \tag{2} $$
where $\Det_{\bX}$~is the {\it determinant line\/} of the Dirac
operator~$D_{\bX}$ on the boundary:
  $$ \Det\mstrut _{\bX} = \bigl(\Det\Ker^-D_{\bX} \bigr)\otimes
     \bigl(\Det\Ker^+D_{\bX} \bigr)\inv . \tag{3} $$
(Recall that $\Det V={\tsize\bigwedge} ^nV$ for an $n$~dimensional vector
space~$V$.  Also $L\inv =L^*$ for a one dimensional vector space~$L$.)
Properly normalized we have $|\tau \mstrut _X|=1$ in the {\it Quillen
metric\/} on~$\dib$.

Now suppose $X\to Z$ is a family of odd dimensional manifolds with boundary.
Then $\bX\to Z$ is a family of closed even dimensional manifolds.  The
determinant lines~\thetag{3} patch together to form a smooth determinant line
bundle $\Det_{\bX/Z}\to Z$.  Furthermore, it carries the Quillen metric and a
canonical connection~$\nabla $, as defined in~\cite{BF1}.  The exponentiated
$\eta $-invariant is now a smooth section
  $$ \tau \mstrut _{X/Z}\:Z\longrightarrow \dibf.  $$

There are two basic results about this invariant: a variation formula and a
gluing law.  The variation formula computes the derivative of~$\tau \mstrut
_{X/Z}$ in a family.

        \proclaim{\protag{4~\cite{DF1,Theorem 1.9}} {Theorem}}
 With respect to the canonical connection~$\nabla $ on~$\dibf$,
  $$ \nabla \tau \mstrut _{X/Z} = \tpi \comp{\int_{X/Z} \Ahat(\Omega
     ^{X/Z})}{1}\cdot \tau \mstrut _{X/Z}.   $$
        \endproclaim

\flushpar
 Here $\Omega ^{X/Z}$~is the Riemannian curvature of~$X\to Z$ and $\Ahat$~is
the usual $\Ahat$-polynomial.  (For other Dirac operators substitute the
appropriate index polynomial in place of~$\Ahat$.)  The `$(1)$'~denotes the
1-form piece of the differential form.  For a family of {\it closed\/}
manifolds this is a result of Atiyah-Patodi-Singer.  The new point here is
the relationship of~$\tau $ with the canonical connection~$\nabla $.  This
plays a crucial role in the next section.

 \midinsert
 \bigskip
 \centerline{
  \epsfxsize= \hsize
 \epsffile{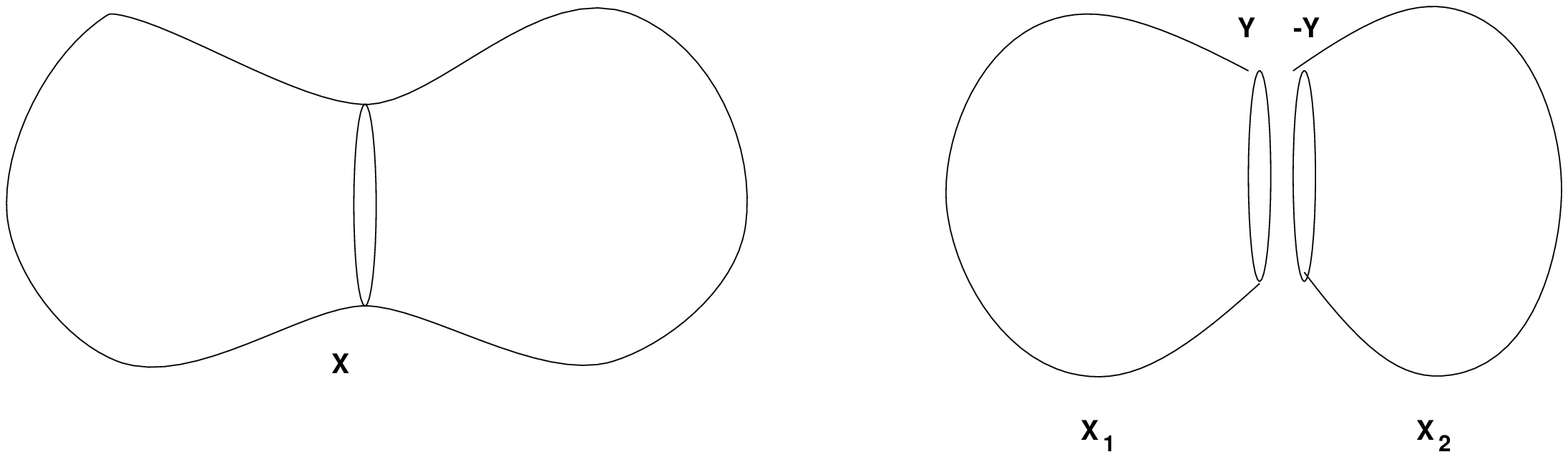}}
 \nobreak
 \botcaption{Figure~1} Cutting a closed manifold into two pieces.
 \endcaption
 \bigskip
 \endinsert

The simplest case of the gluing law is for a {\it closed\/} manifold~$X$
split into two pieces~$X_1,X_2$ along a closed oriented codimension one
submanifold $Y\hookrightarrow X$.  (See Figure~1.)  Then $\tau \mstrut
_{X_i}\in \Det\inv _Y$ and $\tau \mstrut _X\in \CC$.

        \proclaim{\protag{5~\cite{DF1,Theorem 2.20}} {Theorem}}
 In this situation
  $$ \tau \mstrut _X = (\tau \mstrut _{X_1},\tau \mstrut _{X_2})\mstrut
     _{\Det\inv _Y}.  $$
        \endproclaim

 \midinsert
 \bigskip
 \centerline{
 \epsffile{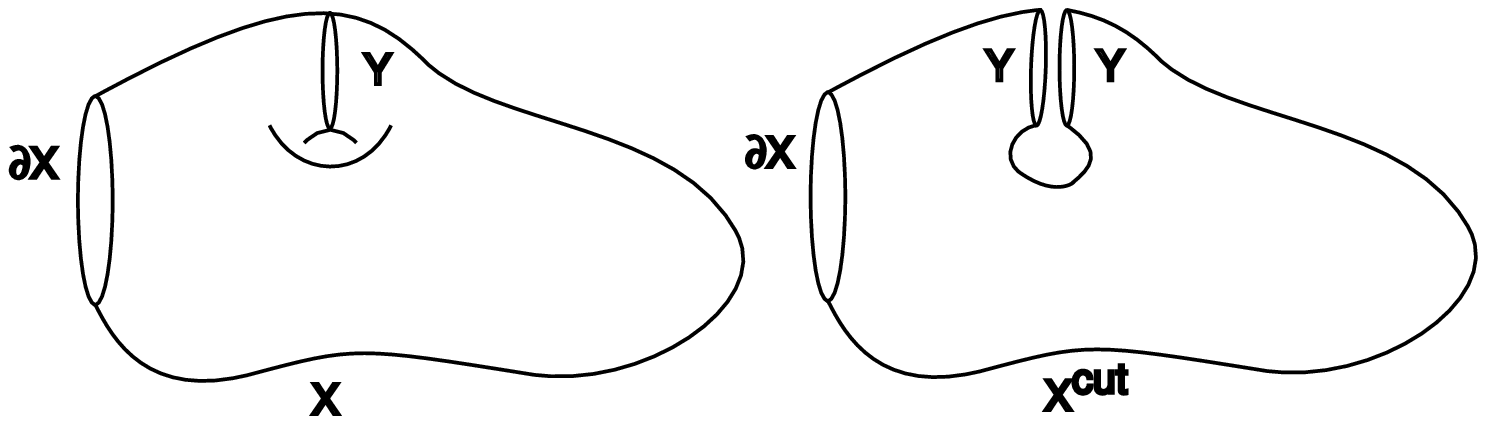}}
 \nobreak
 \botcaption{Figure~2} Cutting a manifold along a submanifold.
 \endcaption
 \bigskip
 \endinsert

The more general gluing formula, which we need in the next section, applies
when $X$~has boundary.  Then for $Y\hookrightarrow X$ a closed oriented
codimension one submanifold we cut along~$Y$ to obtain a new manifold~$X\cut$
with $\bX\cut=\bX\sqcup Y\sqcup -Y$.  (See Figure~2.)  Now
  $$ \aligned
      \tau \mstrut _X&\in \dib\\
      \tau \mstrut _{X\cut}&\in \dib\;\otimes \Det\inv _Y\otimes \Det\inv
     _{-Y}\\
      &\cong \dib\;\otimes \;\;L\mstrut _Y\;\;\;\otimes \;\; L_Y\inv
     ,\endaligned \tag{6} $$
where $L\mstrut _Y=\Det\inv _Y$.  There is now a sign which enters the gluing
formula, and it is nicely taken care of by the following device.  In general
we view the determinant line~$\Det V$ of a vector space~$V$ as a one
dimensional {\it graded\/} vector space whose grading is given by~$\dim V$.
Applied to~\thetag{3} we see that $\Det_Y$ (and so also~$\Det_Y\inv $) is
graded by the {\it index\/} of the Dirac operator~$D_Y$.  Notice that in our
current situation $Y$~does not necessarily bound a 3-manifold, and so its
index may be nonzero.  Let
  $$ \Tr_s\:L\mstrut _Y\otimes L_Y\inv \longrightarrow \CC \tag{7} $$
be the usual contraction times the grading; i.e., if $\index D_Y$ is even it
is the usual contraction and if $\index D_Y$~is odd it is minus the usual
contraction.  That understood, we state the general gluing formula.

        \proclaim{\protag{8~\cite{DF1, Theorem 2.20}} {Theorem}}
 In this situation
  $$ \tau \mstrut _X = \Tr_s(\tau \mstrut _{X\cut}). \tag{9} $$
        \endproclaim

\flushpar
 One of the novel points of~\cite{DF1} is the proof of the gluing law, which
we do not discuss here.

 \subhead Determinant Line Bundles and Adiabatic Limits
 \endsubhead

The application we discuss is to the geometry of the determinant line bundle.
Suppose $\pi \:Y\to Z$ is a family of {\it closed\/} even dimensional
manifolds.  Let $L=\Det\inv _{Y/Z}$ be the inverse determinant line bundle of
the family.  The results in the last section use the Quillen metric and the
construction of the canonical connection~$\nabla $.  But they do not depend
on the formulas for the curvature and holonomy of~$\nabla $, which were
proved in~\cite{BF1}, ~\cite{BF2}.  Here we derive the curvature and holonomy
formulas from \theprotag{4} {Theorem} and \theprotag{8}
{Theorem}.\footnote{As was mentioned in the introduction, this was done in
{}~\cite{DF1,\S5} in an unnecessarily complicated way.  Also, there we {\it
used\/} the curvature formula instead of {\it proving\/} it.  This section
should be considered a rewrite of~\cite{DF1,\S5}.} The basic idea is to use
the $\tau $-invariant~\thetag{2} to define the parallel transport of a new
connection~$\nabla '$ on~$L$.  Thus suppose $\gamma \:\zo\to Z$ is a smooth
path\footnote{Since we need a cylindrical metric near the boundary of~$\Yg$
defined below, we require that $\gamma \bigl([0,\delta ] \bigr)$ and $\gamma
\bigl([1-\delta ,1] \bigr)$ be constant for some~$\delta $.} in~$Z$.  Denote
$I=\zo$.  Let $Y_\gamma =\gamma ^*(Z)\to I$ be the pullback of the family
$\pi \:Y\to Z$ by the path~$\gamma $.  Then $Y_\gamma $~is an odd dimensional
manifold with $\partial Y_\gamma =Z_{\gamma (1)}\sqcup -Z_{\gamma (0)}$.  The
standard metric~$g\mstrut _I$ on~$I=\zo$ determines a metric on~$Y_\gamma $,
since we already have a metric~$g\mstrut _{Y_\gamma /I}$ on the fibers and a
distribution of horizontal planes.  (The projection $\pi \:Y_\gamma \to I$ is
then a Riemannian submersion.)  The $\tau $-invariant of~$Y_\gamma $ is a
linear map
  $$ \tau \mstrut _{Y_\gamma }\:L_{\gamma (0)}\longrightarrow L_{\gamma (1)},
     \tag{10} $$
exactly what we need to define parallel transport.  However, \thetag{10}~does
{\it not\/} define parallel transport since it is not independent of the
parametrization of the path~$\gamma $.  To get a quantity independent of
parametrization we introduce the {\it adiabatic limit\/} as follows.  For
each~$\epsilon \not= 0$ consider the metric
  $$ g\mstrut _\epsilon = \frac{g\mstrut _{I}}{\epsilon ^2}\oplus g\mstrut
_{Y_\gamma /I} \tag{11} $$
on~$Y_\gamma $ relative to the decomposition $T\Yg\cong \pi ^*TI\oplus
T(\Yg/I)$.  Let $\tau \mstrut _{Y_\gamma }(\epsilon )$ be the $\tau
$-invariant for this metric.

        \proclaim{\protag{12} {Lemma}}
 The adiabatic limit
  $$ \tau \mstrut _\gamma =\alim \tau \mstrut _{Y_\gamma } =
     \lim\limits_{\epsilon \to0} \tau \mstrut _{Y_\gamma }(\epsilon )
     \tag{13} $$
exists and is invariant under reparametrization of~$\gamma $.
        \endproclaim

\flushpar
 Notice that the adiabatic limit is introduced for a simple geometrical
reason---to scale out the dependence of~$\tau $ on the parametrization
of~$\gamma $.

        \demo{Proof}
 Here we follow~\cite{DF1,\S5}.\footnote{And we correct a mistake in the
exposition there.} As a preliminary we state without proof a simple result
about the Riemannian geometry of adiabatic limits.  Let $\nabla ^{Y_\gamma
}(\epsilon )$~denote the Levi-Civita connection on~$\Yg$ of the
metric~\thetag{11} and $\curv{\Yg}(\epsilon )$~its curvature.  The result we
need, which follows from a straightforward computation in local Riemannian
geometry, is that $\alim\nabla ^{\Yg}=\lim\limits_{\epsilon \to0}\nabla
^{\Yg}(\epsilon )$~exists and is torsionfree.  Furthermore, the curvature of
this limiting connection is the limit of the curvatures of~$\nabla
^{\Yg}(\epsilon )$ and has the form
  $$ \alim\curv{\Yg}=\lim\limits_{\epsilon \to0}\curv{\Yg}(\epsilon ) =
     \pmatrix 0&*\\ 0&\curv{\Yg/I} \endpmatrix  $$
relative to a fixed (nonorthonormal) basis.  It follows that
  $$ \alim\Ahat(\curv{\Yg}) = \lim\limits_{\epsilon \to0}\Ahat(\curv{\Yg}) =
     \Ahat(\curv{\Yg/I}). \tag{14} $$
We apply this result to families of adiabatic limits, where it also holds.

To prove that the adiabatic limit exists, consider the family of Riemannian
manifolds $\Yg\times \Rnz\to\Rnz$, where the metric on the fiber at~$\epsilon
$ is~\thetag{11}.  According to the variation formula \theprotag{4} {Theorem}
we have
  $$ \frac{d}{d\epsilon }\tau \mstrut _{\Yg}(\epsilon ) =
     \tpi\comp{\int_{(\Yg\times \Rnz)/\Rnz}\Ahat\bigl(\curv {(\Yg\times
     \Rnz)/\Rnz}\bigr)}{1}. \tag{15} $$
Now \thetag{14}~implies that
  $$ \lim\limits_{\epsilon \to0} \Ahat\bigl(\curv {(\Yg\times
     \Rnz)/\Rnz}\bigr) = \Ahat(\curv{\Yg/I}). \tag{16} $$
One should understand this as a limit of sections of a bundle on~$\Yg$ whose
fibers are forms on $\Yg\times \{0\}$.  In other words, they are forms
on~$\Yg$ with a `$d\epsilon $' term as well.  Formula~\thetag{16} implies
that there is no $d\epsilon $~term in the limit, and so the integral over the
fibers in~\thetag{15} vanishes.  Therefore, $\lim\limits_{\epsilon
\to0}\frac{d}{d\epsilon }\tau \mstrut _{\Yg}(\epsilon )=0$ and so $\alim\tau
\mstrut _{\Yg}=\lim\limits_{\epsilon \to0}\tau \mstrut _{\Yg}(\epsilon )$
exists.

A similar argument proves that $\tau \mstrut _\gamma $ is invariant under
reparametrization.  Let $\Cal{D}$~denote the space of diffeomorphisms $\phi
\:\zo\to\zo$ with $\phi (0)=0$ and~$\phi (1)=1$.  We pull back $\pi \:Y\to Z$
via the map
  $$ \aligned
      \zo\times \Rnz\times \Cal{D}&\longrightarrow Z\\
      \langle t,\epsilon ,\phi \rangle&\longmapsto \gamma \bigl(\phi (t)
     \bigr)\endaligned  $$
to construct the family of manifolds
  $$ \Cal{Y}\longrightarrow \Rnz\times \Cal{D}, $$
where the metric on the fiber over~$\langle \epsilon ,\phi \rangle$
is~\thetag{11}.  As in the previous argument we compute the differential
of~$\tau \mstrut _{Y_{\gamma \circ \phi }}(\epsilon ,\phi )$ in the adiabatic
limit:
  $$ \lim\limits_{\epsilon \to0}d\tau \mstrut _{\langle \epsilon ,\phi
     \rangle} = \tpi\;\;\sigma ^*\!\comp{\int_{Y/Z}\Ahat(\curv{Y/Z})}{2},
     \tag{17} $$
where
  $$ \aligned
      \sigma \:\zo\times \Cal{D}&\longrightarrow Z\\
      \langle t,\phi \rangle&\longmapsto \gamma \bigl(\phi (t) \bigr)
     \endaligned $$
We conclude that \thetag{17}~vanishes since the image of~$\sigma $ is one
dimensional---the pullback of a 2-form vanishes.
        \enddemo

        \proclaim{\protag{18} {Lemma}}
 The maps~$\tau \mstrut _\gamma $ are the parallel transport of a
connection~$\nabla '$ on $L\to Z$.
        \endproclaim

        \remark{Remark}
 Since $\tau \mstrut _\gamma $~is a unitary transformation ($|\tau \mstrut
_\gamma |=1$), the connection~$\nabla '$ is also unitary.
        \endremark

        \demo{Proof}
 By a general result~\cite{F2, Appendix~B} it suffices to show that the
fiducial parallel transport $\tau \mstrut _\gamma $~is invariant under
reparametrization and composes under gluing.  The first statement is
contained in the previous lemma.  For the second, if $\gamma _1,\gamma _2$
are paths with $\gamma _2(0) = \gamma _1(1)$, then we can compose to get a
path~$\gamma =\gamma _2\circ \gamma _1$.  The gluing law \theprotag{8}
{Theorem} then implies $\tau \mstrut _\gamma =\tau \mstrut _{\gamma _2}\circ
\tau \mstrut _{\gamma _1}$ as required.  (\theprotag{8} {Theorem} applies to
a fixed metric and then we take the adiabatic limit.)
        \enddemo

        \remark{Remark}
 It is instructive to see in detail how the sign works in this application of
the gluing law.  Here we cut~$\Yg$ along $Y=Y_{\gamma _2(0)}= Y_{\gamma
_{1}(1)}$ to obtain~$\Yg\cut=Y_{\gamma _1}\sqcup Y_{\gamma _2}$.  So
  $$ \aligned
      \tau \mstrut _{\gamma _1}&\in \Hom(L_{\gamma _1(0)}, L_{\gamma _1(1)})
     \cong \;\;L\mstrut _Y\;\;\otimes L\inv _{\gamma _1(0)}, \\
      \tau \mstrut _{\gamma _2}&\in \Hom(L_{\gamma _2(0)}, L_{\gamma _2(1)})
     \cong L _{\gamma _2(1)}\otimes L _Y\inv , \endaligned  $$
where we write $L_Y = L_{\gamma _1(1)} = L_{\gamma _2(0)}$.  Thus
  $$ \tau \mstrut _{\Yg\cut} = \tau \mstrut _{\gamma _2}\otimes \tau \mstrut
     _{\gamma _1} \in L_{\gamma _2(1)}\otimes L_Y\inv \otimes L\mstrut
     _Y\otimes L\inv _{\gamma _1(0)}. \tag{19} $$
The key point is that the factors are in a different order than
in~\thetag{6} and~\thetag{7}---now the factor~$L_Y\inv $ {\it precedes\/} the
factor~$L\mstrut _Y$.  So the contraction is the usual trace.  Put
differently, to move~\thetag{19} to the standard form~\thetag{6} we
introduce a factor of~$\sind$ and this is cancelled by the factor~$\sind$ in
the supertrace~\thetag{9}.  The upshot is that in this situation the right
hand side of~\thetag{9} is~$\tau \mstrut _{\gamma _2}\circ \tau \mstrut
_{\gamma _1}$ as desired.
        \endremark

It is quite easy to prove from the variation formula \theprotag{4} {Theorem}
that this new connection agrees with the canonical connection~$\nabla $.

        \proclaim{\protag{20} {Proposition}}
 $\nabla '=\nabla $.
        \endproclaim

        \demo{Proof}
 We must show that the parallel transports agree.  Let $\gamma \:\zo\to Z$ be
a path and fix an element $\ell _0\in L_{\gamma (0)}$ of unit norm.  Then if
$\gamma \:[0,t]\to Z,\ 0\le t\le1,$ is the restriction of~$\gamma $, and
$\tau \mstrut _t\:L_{\gamma (0)}\to L_{\gamma (t)}$ the parallel transport
of~$\nabla '$, by definition the path $\ell _t=\tau \mstrut _t(\ell _0)$ is
parallel for~$\nabla '$.  It suffices to show that $\dfrac{D\tau \mstrut
_t}{Dt}=0$, where $\dfrac{D}{Dt}=\nabla $ along the path~$\gamma $.  For then
$\dfrac{D\tau \mstrut _t(\ell _0)}{Dt}=0$ as well, since $\ell _0$~is a
constant.

Define $T=\{\langle t,s  \rangle\in \zo\times \zo:s\le t\}$ with projection
  $$ \aligned
       \rho \:T&\longrightarrow \zo=I\\
       \langle t,s \rangle&\longmapsto t\endaligned  $$
and a map
  $$ \aligned
      \Gamma \:T&\longrightarrow Z\\
      \langle t,s  \rangle&\longmapsto \gamma (s).  \endaligned  $$
Then the pullback $\pi \:\Gamma ^*Y\to T$ determines a family of manifolds
$\rho \circ \pi \:\Gamma ^*Y\to\zo$ parametrized by~$I=\zo$.  We use the flat
metric on~$T$ and make $\pi \:\Gamma ^*Y\to T$ a Riemannian submersion.  The
variation formula \theprotag{4} {Theorem} implies
  $$ \frac{D\tau \mstrut _t}{Dt} = \tpi \int_{\Gamma ^*Y/I}
     \alim\comp{\Ahat(\Omega ^{\Gamma ^*Y/I})}{1}. \tag{21} $$
Even before taking the adiabatic limit, the fact that $\Gamma $~factors
through the projection $\langle t,s  \rangle\mapsto s$ implies that the right
hand side of~\thetag{21} vanishes.
        \enddemo

In view of \theprotag{20} {Proposition}, to compute the curvature and holonomy
of~$\nabla $ it suffices to compute the curvature and holonomy of~$\nabla '$.
Notice that since $L=\Det\inv _{Y/Z}$~is the {\it inverse\/} determinant line
bundle our formulas here have opposite signs to those for~$\Det_{Y/Z}$
computed in~\cite{BF1}, ~\cite{BF2}.  The
holonomy is computed from the parallel transport by a straightforward
application of the gluing law.  We must only be careful about the spin
structure.  Recall that $\cir$~has two spin structures.  The {\it
nonbounding\/} spin structure is the trivial double cover of the circle; the
{\it bounding\/} spin structure is the nontrivial double cover.

        \proclaim{\protag{22~\cite{BF2,Theorem 3.18}} {Theorem}}
 Suppose $\gamma \:\zo\to Z$ is a {\it closed\/} path.\footnote{Recall that we
require that $\gamma \bigl([0,\delta ] \bigr)$ and $\gamma \bigl([1-\delta
,1] \bigr)$ be constant for some~$\delta $.} There is an induced
manifold~$\Ygh\to\cir$ obtained by gluing the ends of~$\Yg$.  Then the
holonomy of~$L$ around~$\gamma $ is
  $$ \hol_L(\gamma ) =\cases \sind\alim\tau \mstrut _{\Ygh}
     ,&\text{nonbounding spin structure on $\cir$};\\\alim\tau \mstrut
     _{\Ygh},&\text{bounding spin structure on $\cir$}.\endcases \tag{23} $$
        \endproclaim

\flushpar
 Here the spin structure on~$\cir$ combines with the spin structure
on~$T(\Ygh/\cir)$ to give a spin structure on~$\Ygh$.

        \demo{Proof}
 This follows directly from the definition~\thetag{13} of parallel transport
and the gluing law applied to~$X=\Ygh$ and~$X\cut=\Yg$.  Take first the
nonbounding spin structure on~$\cir$, lifted to a spin structure on~$\Ygh$.
The induced spin structure on the cut manifold~$\Yg$ is the standard one,
with the ends each identified with $Y_z$, where $z=\gamma (0)=\gamma (1)$.
Now for each~$\epsilon $ the $\tau $-invariant of~$\Yg$ is an element
  $$ \tau \mstrut _{\Yg(\epsilon )}\in L\mstrut _z\otimes L_z\inv .
      $$
Then \theprotag{8} {Theorem} implies
  $$ \tau \mstrut _{\Ygh(\epsilon )}=\sind\tau \mstrut _{\Yg(\epsilon )},
      $$
where on the right hand side we identify $L\mstrut _z\otimes L_z\inv $
with~$\CC$ using the {\it usual contraction\/}.  Now the first equation
in~\thetag{23} follows from the definition of holonomy in terms of parallel
transport.  To obtain the second equation, consider the identity map of~$Y_z$
lifted to the {\it nontrivial\/} deck transformation on the spin bundle
of~$Y_z$.  It induces multiplication by~$\sind$ on the inverse determinant
line~$L_z$.  Apply this transformation to~$\Yg$ before gluing in order to
switch spin structures on~$\Ygh$.  Then the second equation in~\thetag{23}
follows from the first.
        \enddemo

        \proclaim{\protag{24~\cite{BF2,Theorem~1.21}} {Theorem}}
 The curvature~$\curv L$ of the inverse determinant line bundle $L\to Z$ is
  $$ \curv L = -\tpi\comp{\int_{Y/Z}\Ahat(\curv{Y/Z})}{2}. \tag{25} $$
        \endproclaim

        \demo{Proof}
 For any line bundle we can determine the curvature once we know the holonomy
as follows.  Suppose $\Gamma \:D\to Z$ is a map of a disk into~$Z$ with
boundary map~$\gamma $.  Let $\YG=\Gamma ^*Y\to D$ be the pullback manifold;
then $\partial \YG = \Ygh$.  In the following calculation we use the bounding
spin structure on~$\cir$ and the induced spin structure on~$\Ygh$.
  $$ \aligned
      \int_{D}\curv L &= -\log\hol_L(\gamma ),\\
      &= \alim\bigl(-\log \tau \mstrut _{\Ygh} \bigr),\\
      &= \alim\left\{ -\tpi\int_{\YG}\Ahat(\curv{\YG})  \right\},\\
      &=\int_{D}(-\tpi)\int_{\YG/D}\Ahat(\curv{\YG/D}),\\
      &=\int_{D}\Gamma ^*\left\{ -\tpi\comp{\int_{Y/Z}\Ahat(\curv{Y/Z})}{2}
     \right\} .\endaligned \tag{26} $$
In the fourth line we apply~\thetag{14}.  In the third line we apply the
index theorem of Atiyah-Patodi-Singer~\cite{APS} which asserts that
  $$ \int_{\YG}\Ahat(\curv{\YG}) - \frac{\eta \mstrut _{\YG}(0)+\dim\Ker
     D_{\YG}}{2}  $$
is a certain index, so in particular is an integer.  When $\Gamma $~shrinks
the disk to a point both sides of~\thetag{26} vanish, so we have chosen the
correct logarithm on the right hand side of~\thetag{26}.  Since~\thetag{26}
holds for all $\Gamma \:D\to Z$, equation~\thetag{25} follows.
        \enddemo


\Refs\tenpoint

\ref
\key APS
\by M. F. Atiyah, V. K. Patodi, I. M. Singer
\paper Spectral asymmetry and Riemannian geometry. I
\jour Math. Proc. Cambridge Philos. Soc. \vol 77 \yr 1975 \pages 43--69
\endref

\ref
\key BF1
\by J. M. Bismut, D. S. Freed \paper The analysis of elliptic
families I: Metrics and connections on determinant bundles \jour Commun. Math.
Phys. \vol 106 \pages 159--176 \yr 1986
\endref

\ref
\key BF2
\by J. M. Bismut, D. S. Freed \paper The analysis of elliptic families II:
Dirac operators, eta invariants, and the holonomy theorem of Witten \jour
Commun. Math. Phys. \vol 107 \yr 1986 \pages 103--163
\endref

\ref
\key DF1
\by X. Dai, D. S. Freed
\paper $\eta $-invariants and determinant lines
\jour J. Math. Phys.
\yr 1994
\vol 35
\pages 5155--5194
\endref

\ref
\key DF2
\by X. Dai, D. S. Freed
\paper $\eta $-invariants and determinant lines
\jour C. R. Acad. Sci. Paris
\yr 1995
\pages 585--592
\endref

\ref
\key F1
\by D. S. Freed \paper On determinant line bundles \inbook Mathematical Aspects
of String Theory \bookinfo ed. S. T. Yau \publ World Scientific Publishing \yr
1987
\endref

\ref
\key F2
\by D. S. Freed
\paper Classical Chern-Simons theory, Part 1
\jour Adv. Math.
\toappear
\endref

\ref
\key Q
\by D. Quillen \paper Determinants of Cauchy-Riemann operators over a Riemann
surface \jour Funk. Anal. iprilozen \vol 19 \yr 1985 \pages 37
\endref

\ref
\key W
\by E. Witten \paper Global gravitational anomalies \jour Commun. Math. Phys.
\vol 100 \yr 1985 \pages 197--229
\endref

\endRefs
\enddocument